\documentclass[12pt]{amsart}
\usepackage{amsfonts}
\usepackage{amssymb}
\theoremstyle{definition}

\theoremstyle{remark}

\newcommand{\ds}{\displaystyle}

\begin{document}

\centerline{Riv. Mat. Univ. Parma (5) {\bf 4} (1995), 161-168}

\vspace{0.5in}
\centerline{{\large O}GNIAN  {\large K}ASSABOV
\footnote{Higher Transport School "T. Kableschkov", Section of Mathematics, Slatina, 1574 Sofia, Bulgaria.}}

\vspace{0.3in}

\centerline{\bf On conformally flat totally real submanifolds
\footnote{Received June 1, 1995. AMS classification 53C40. The research has bees partially
supported by the Ministry of Education of Bulgaria, Contract MM 413/94.}}

\vspace{0.6in}
{\bf 1 - Introduction}

\vspace{0.2in}
Let $\widetilde {M}^{2m}$ be a $2m$-{\it dimensional K\"ahler manifold} \, with Riemannian metric g,
complex structure $J$ and Riemannian connection $\widetilde \nabla$. The curvature tensor, the Ricci 
tensor and the scalar curvature of $\widetilde {M}^{2m}$ are denoted by $\widetilde R$, $\widetilde S$,
$\tilde \tau$, respectively. The Bochner curvature tensor $\widetilde B$ of $\widetilde M $ is given by 

$$
	\widetilde B = \widetilde R - \frac 1{2(m+2)}(\varphi + \psi)(\widetilde S) +
	\frac{\tilde\tau}{8(m+1)(m+2)}(\varphi + \psi) (g)
$$

\noindent
where the operators $\varphi$ and $\psi$ are defined by

$$
	\begin{array}{rl}
	\varphi (Q)(x,y,z,u)& = g(x,u)Q(y,z) - g(x,z)Q(y,u)
	                     + g(y,z)Q(x,u) - g(y,u)Q(x,z)  \\
	                     & \\
	\psi (Q)(x,y,z,u) &	= g(x,Ju)Q(y,Jz) - g(x,Jz)Q(y,Ju) - 2g(x,Jy)Q(z,Ju)  \\
	                  & +\, g(y,Jz)Q(x,Ju) - g(y,Ju)Q(x,Jz) - 2g(z,Ju)Q(x,Jy)
	 \end{array}
$$

\noindent
for a tensor field $Q$ of type (0,2), and $x,\, y,\, z,\, u$ are vector fields of  
$\widetilde {M}^{2m}$.

Let $M$ be a {\it submanifold of} $\widetilde {M}$. The Gauss and Weingarten formulas are
given respectively by

$$
	\widetilde\nabla_X Y = \nabla_X Y +\sigma(X,Y)  \qquad 
	\widetilde\nabla_X \xi = -A_\xi X + D_X \xi
$$  

\noindent
for vector fields \ $X,\, Y$ \ tangent to \ $M$ \ and \ $\xi$ \ normal to \ $M$, where \ $\nabla$ \ is
the Riemannian connection on \ $M$, \ $D$ is the normal connection, \ $\sigma$ \ is the second
fundamental form of \ $M$ \ and \ $A_\xi X$ \ is the tangential component of \ $\widetilde\nabla_X \xi$. 
It is well known that \
$ g(\sigma(X,Y),\xi) = g(A_\xi X,Y)$. The mean curvature vector $H$ is defined by
$H=\ds\frac{1}{n} {\rm tr} \, \sigma$. If $H=0$, $M$ is called minimal. In particular,
if $\sigma = 0$, $M$ is said to be  totally geodesic. A normal vector field $\xi$ is said
to be parallel, if $ D_X \xi = 0$ for each vector field $X$ on $M$. 

An $n$-dimensional submanifold $M^n$ of  $\widetilde {M}^{2m}$ is said to be a {\it totally
real submanifold} of  $\widetilde {M}^{2m}$, if for each point $p \in M^n$,
$ JT_pM^n \subset T_p(M^n)^{\perp}$. Then $n \le m$. In the following we suppose that $M^n$
{\it is a totally real submanifold of}  $\widetilde {M}^{2m}$  and $m=n$. In this case it is 
not difficult to find 

$$
	\sigma (X,Y) = JA_{JX} Y = JA_{JY} X      \leqno (1.1)
$$
$$
	D_X JY = J\nabla_X Y \, .                \leqno (1.2)
$$

\noindent
See e.g. \cite{YK2}.

Let $\widetilde R$ be the curvature tensor of  $\widetilde {M}^{2n}$. Then using (1.1),
the Gauss equation can be written as

$$
	\{ \widetilde R(X,Y)Z \}^t = R(X,Y)Z - [A_{JX},A_{JY}]Z
$$

\noindent
where $t$ denotes the tangential component.

Let $\overline \nabla$ denote the connection of van der Waerden-Bortolotti.
Then $M$ is said to be a {\it parallel submanifold} of $\widetilde M$ if
$\overline \nabla\sigma = 0$. More generally $M$ is called a {\it semiparallel
submanifold} of $\widetilde M$, if \ $\overline R(X,Y).\sigma=0$, where

$$
	(\overline R(X,Y).\sigma)(Z,U) = R^{\perp}(X,Y)\sigma(Z,U) -
	   \sigma(R(X,Y)Z,U) - \sigma(Z,R(X,Y)U) 
$$

\noindent
$R^\perp$ being the curvature tensor of the normal connection $D$. The
investigation of semiparallel submanifolds was initiated by J. Deprez \cite{D}.
For a semiparallel submanifold by (1.1) and (1.2) we obtain

$$
	R(X,Y)A_{JZ}U = A_{JZ} R(X,Y)U + A_{JU} R(X,Y)Z \, .    \leqno (1.3)
$$

\noindent
On the other hand the {\it submanifolds with semiparallel mean curvature
vector} are defined by $ R^\perp (X,Y)H = 0$ \cite{OK}. Note that the class
of submanifolds with semiparallel mean curvature includes the semiparallel submanifolds and the
submanifolds with parallel mean curvature vector. 

Let $S$ and $\tau$ denote the Ricci tensor and the scalar 
curvature of $M$, respectively. Then as it is well known for $n > 3$ \ $M$ is 
{\it conformally flat}, if and only if the Weil conformal curvature tensor $C$ of 
$M$ vanishes, where

$$
	C =  R - \frac 1{n-2}\varphi (S) +
	\frac{\tau}{2(n-1)(n-2)}\varphi (g) \, .
$$

In section 2 and 3 we prove

\vspace{0.1in}
{\ttfamily Theorem 1.} {\it Let $M^n$ be a conformally flat totally real submanifold of a 
K\"ahler manifold  $\widetilde {M}^{2n}$, $n>3$. Assume also that the mean curvature vector of
$M^n$ is semiparallel. If $M^n$ is not minimal at a point $p$, then, in a neighborhood 
of $p$, $M^n$ is either flat or a product $M_1^{(n-1)}(c) \times I$, where $M_1^{(n-1)}(c)$ is an 
$(n-1)$-dimensional manifold of constant sectional curvature $c \ne 0$ and $I$ is a segment.}

\vspace{0.1in} 
{\ttfamily Theorem 2.} {\it Let $M^n$ be a conformally flat totally real semiparallel submanifold of a 
K\"ahler manifold  $\widetilde {M}^{2n}$, $n>3$. If $M^n$ is not totally geodesic at a
point $p$, then, in a neighborhood of $p$, $M^n$ is flat or $M^n=M_1^{(n-1)}(c) \times I$.}

\vspace{0.1in}
In section 4 we deal with products of K\"ahler manifolds with vanishing 
Bochner curvature tensor.

\vspace{0.4in}
{\bf 2 - Proof of Theorem 1}

\vspace{0.2in}

First we prove

\vspace{0.1in}
{\ttfamily Proposition.} {\it Let $M^n$ be a conformally flat totally real submanifold of a K\"ahler
manifold $\widetilde {M}^{2n}$, $n>3$, such that the mean curvature vector $ H $ is
semiparallel at a point $p$. If $M^n$ is not minimal at $p$, it is quasi-Einstein at $p$
with $ S(JH_p,JH_p)=0$.} 

\vspace{0.1in}
{\ttfamily Proof.} Let $\{ e_i \} \ i=1,...,n$ be an orthonormal basis of $T_pM$, such that 
$ S_{1,1}(e_i) = \lambda_ie_i $ for $ i = 1,...,n$, where $ S_{1,1} $ is the Ricci
tensor of type (1,1). From $ C=0 $ and  $ R(e_i,e_j)JH=0 $ we obtain 

$$
	(\lambda_i+\lambda_j - \frac{\tau}{n-1})g(e_j,JH) = 0 \ .   \leqno (2.1)
$$

\noindent
If $ g(e_j,JH) = 0 $ for each $j$, $M^n$ is minimal at $p$. Let e.g. $g(e_1, JH) \ne 0$.
Then (2.1) implies

$$
	\lambda_1+\lambda_i - \frac{\tau}{n-1} = 0         \leqno (2.2)
$$

\noindent
for $ i=2,...,n$. Hence $\lambda_i = \lambda_j $ for $i,j=2,...,n$, i.e. $M^n$ is
{\it quasi-Einstein} at $p$. Moreover (2.2) implies $ \lambda_1 = 0 $. If there exists
$ i>1 $, such that \ $g(e_i,JH) \ne 0$, it follows that $M^n$ is Einstein at $p$, with
$\tau =0$, so $S=0$ at $p$. If \ $g(e_i,JH) = 0$ \ for $i=2,...,n$ \ it follows that $(JH)_p$
is proportional to $e_1$, thus proving our assertion.

Now we can prove Theorem 1. Since $H$ does not vanish at $p$, then this holds
also in a neighborhood of $p$. Then Theorem 1 follows from our Proposition and a theorem 
of Kurita, see \cite{K}.

If $H$ has constant lenght, Then $M$ is minimal or $H$ does not vanishes. Hence we have:

\vspace{0.1in}
{\ttfamily Corollary 1.} {\it Let $M^n$ be a conformally flat totally real submanifold of a K\"ahler 
manifold $\widetilde {M}^{2n}$. Assume also that the mean curvature vector $H$ of $M$
is semiparallel and with constant lenght. Then one of the following holds:

$M^n$ is minimal

$M^n$ is locally flat or a product $M_1^{n-1}(c) \times I, \ c \ne 0$.

In particular the result is true when $H$ is parallel.}

\vspace{0.4in}
{\bf 3 - Proof of Theorem 2}

\vspace{0.2in}
As in Section 2 \ let \ $\{ e_i \} \ i=1,...,n$ \ be an orthonormal basis of \ $T_pM$, such that 
$ S_{1,1}(e_i) = \lambda_ie_i $ \ for \ $ i = 1,...,n$. Under the assumptions of Theorem 2 
we prove some lemmas.

\vspace{0.1in}
{\ttfamily Lemma 1.} {\it Let there exist \ $i \ne k$, such that \ $g(A_{Je_i} e_i,e_k) \ne 0$. Then
$\lambda_i=\lambda_j$ for all \ $j \ne k$ \ and \ $\lambda_k = 0$.} 

\vspace{0.1in}
{\ttfamily Proof.} We put in (1.3) \ $X=e_j,\, Y=e_k,\, Z=U=e_i$ \ for \ $j\ne i,k$ \ and we obtain

$$
	(\lambda_j+\lambda_k - \frac{\tau}{n-1})g(A_{Je_i}e_i,e_k) = 0 
$$ 

\noindent
which implies

$$
	\lambda_j+\lambda_k - \frac{\tau}{n-1}  \ . \leqno (3.1)
$$

\noindent
Now we put in (1.3) $X=e_k,\, Y=Z=U=e_i$ and we find

$$
	(\lambda_i+\lambda_k - \frac{\tau}{n-1})(2A_{Je_i}e_k+g(A_{Je_i} e_i,e_k)e_i
	                 -g(A_{Je_i}e_i,e_i)e_k) = 0         \leqno (3.2) 
$$ 

\noindent
which implies

$$
	\lambda_i+\lambda_k - \frac{\tau}{n-1} = 0 \ .        \leqno (3.3) 
$$

\noindent
From (3.1) and (3.3) it follows $\lambda_i=\lambda_j$. Then (3.1) implies $\lambda_k=0$.

\vspace{0.1in}
{\ttfamily Lemma 2.} {\it Let there exists \ $i$, such that \ $g(A_{Je_i} e_i,e_i) \ne 0$. Then $\lambda_i=0$
and \ $\lambda_j=\lambda_k$ \ for all \ $j,k \ne i$.}   

\vspace{0.1in}
{\ttfamily Proof.} If we have

$$
	\lambda_i+\lambda_k - \frac{\tau}{n-1} = 0
$$

\noindent
for any $k=1,...,n$ the assertion follows immediately. Let us assume that there exists a
$k$ such that

$$
	\lambda_i+\lambda_k - \frac{\tau}{n-1} \ne 0 \ .
$$

\noindent
As in Lemma 1 we find (3.2) and hence we have

$$
	2A_{Je_i}e_k+g(A_{Je_i} e_i,e_k)e_i-g(A_{Je_i}e_i,e_i)e_k = 0         
$$ 

\noindent
which implies \
$ g(A_{Je_k}e_k,e_i) \ne 0$. Using Lemma 1 we obtain $\lambda_i=0$ and $\lambda_j=\lambda_k$
for $j,k \ne i$.

\vspace{0.1in}
{\ttfamily Lemma 3.} {\it Let \ $M$ \ be minimal at \ $p$. Then \ $M$ \ is totally geodesic at \ $p$ \ or there exists
\ $k$, such that \ $\lambda_k=0$ \ and $\lambda_i=\lambda_j$ for \ $i,j \ne k$.}

\vspace{0.1in}
{\ttfamily Proof.} If there exists \ $i$, such that \ $A_{Je_i}e_i \ne 0$, the assertion follows from
Lemmas 1 and 2. So let \ $A_{Je_i}e_i = 0$ \ for any \ $i=1,...,n$. Suppose that $M$ is not
totally geodesic at $p$. Then \ $g(A_{Je_i}e_j,e_k) \ne 0$ \ for some \ $i\ne j \ne k \ne i$.
We put in (1.3) \ $X=U=e_s$, $Y=e_j$, $Z=e_k$ \ for \ $s\ne j,k$ \ and we obtain

$$
	(\lambda_s+\lambda_j - \frac{\tau}{n-1})(A_{Je_j}e_k+g(A_{Je_s} e_k,e_j)e_s) = 0   
$$ 

\noindent
which implies

$$
	\lambda_s+\lambda_j - \frac{\tau}{n-1} = 0 \ .   
$$ 

\noindent
Analogously

$$
	\lambda_s+\lambda_k - \frac{\tau}{n-1} = 0    
  \qquad
	\lambda_t+\lambda_k - \frac{\tau}{n-1} = 0    
$$ 

\noindent
for \ $t \ne i,\, k$. Hence it follows \ $\lambda_l=0$ \ for any \ $l=1,\hdots,n$, which
proves the Lemma.

Now we are in position to prove Theorem 2. If \ $M$ \ is not minimal at \ $p$, Theorem 2 
follows from Theorem 1. Let \ $M$ \ be minimal at \ $p$. Then the assertion follows 
from Lemma 3 and \cite{K}.

\vspace{0.4in}
{\bf 4 - Submanifolds of Bochner flat K\"ahler products}

\vspace{0.2in}
Let \ $\widetilde {M}^{2n}$ \ be a  K\"ahler manifold with vanishing Bochner
curvature tensor and constant scalar curvature. Then \ $\widetilde {M}^{2n}$ \
either has constant holomorphic sectional curvature or is locally a product of two 
K\"ahler manifolds of constant holomorphic sectional curvature $\mu$ and $-\mu$,
respectively, $\mu>0$, \cite{MT}. Totally real submanifolds of K\"ahler manifolds
of constant holomorphic sectional curvature have been studied by many authors,
see e.g. \cite{CO}, \cite{YK1}, \cite{YK2}. Now we consider the case of K\"ahler products 
with vanishing Bochner curvature tensor.  

\vspace{0.1in}
{\ttfamily Theorem 3.} {\it Let \, $M^n$ \, be a totally real semiparallel submanifold with commutative
second fundamental form and mean curvature vector of constant lenght of a K\"ahler
product \, $\widetilde {M}^{2k}(\mu) \times \widetilde {M}^{2(n-k)}(-\mu)$, $\mu \ne 0$,
$n>3$, $k\ge n-k\ge 1$. Then $M^n$ is a product \, $M^k(\ds\frac \mu 4) \times M^{n-k}$,
where \, $M^k(\ds\frac \mu 4)$ \, is a manifold of constant curvature \, $\ds \frac{\mu}4$ \, and 
is totally geodesic in \ $\widetilde {M}^{2k}(\mu)$. If in addition \ $n-k>1$, then \,
$M^{n-k}$ \, is totally geodesic in \, $ \widetilde {M}^{2(n-k)}(-\mu)$ \, and has constant
sectional curvature \, $\ds -\frac \mu 4$.}

\vspace{0.1in}
{\ttfamily Proof.} Since $M^n$ has commutative second fundamental form (i.e.
$A_{\xi}A_{\eta} = A_{\eta}A_{\xi}$ $\forall \xi,\, \eta \in TM^\perp$, \cite{YK2}, p.29),
the Gauss equation implies

$$
	\widetilde R(X,Y,Z,U) = R(X,Y,Z,U)
$$

\noindent
for arbitrary vectors \ $X,\, Y,\, Z,\, U$ \ in \ $T_pM$. Let \ $X,\, Y,\, Z,\, U$ \ be orthogonal.
Since $\widetilde B=0$ we obtain $ R(X,Y,Z,U)=0$ and hence $M^n$ is conformaly flat, see
e.g. \cite{S} p. 307. If $M^n$ is totally geodesic, it is straightforward that it is a
product $\ds M^k(\frac{\mu}{4}) \times M^{n-k}(-\frac{\mu}{4})$, where $\ds M^k(\frac{\mu}{4})$,
resp. $\ds M^{n-k}(-\frac{\mu}{4})$, is totally geodesic in $\widetilde {M}^{2k}(\mu)$,
resp. $ \widetilde {M}^{2(n-k)}(-\mu)$.

Let \ $M^n$ \ be not totally geodesic. According to Theorem 2 it is locally flat or a product \
$\ds M_1^{n-1}(c) \times I$. As it is easily seen, if \ $M^n$ is flat, it follows \ $\mu = 0$, 
which is not our case. So \ $M^n$ \ is locally \ $\ds M_1^{n-1}(c) \times I$. Denote by \ $\pi_1$ \
and \ $\pi_2$  the projections of \ $\widetilde {M}^{2k}(\mu) \times \widetilde {M}^{2(n-k)}(-\mu)$,
onto \ $\widetilde {M}^{2k}(\mu) $ \ and \ $ \widetilde {M}^{2(n-k)}(-\mu)$, respectively. 
The induceed differentials will be denoted also by $\pi_1$ and $\pi_2$. Let \ $F=\pi_1-\pi_2$.
Then we have \cite{M}, \cite{T}

$$
	\begin{array}{c}
		\widetilde R(\tilde x,\tilde y,\tilde z,\tilde u) = \ds\frac{\mu}{8} 
		\{ g(F\tilde x,\tilde u)g(\tilde y,\tilde z) - g(F\tilde x,\tilde z)g(\tilde y,\tilde u) \\
		+ g(\tilde x,\tilde u)g(F\tilde y,\tilde z) - g(\tilde x,\tilde z)g(F\tilde y,\tilde u)
							+ g(J\tilde x,\tilde u)g(JF\tilde y,\tilde z) \\   
		- g(J\tilde x,\tilde z)g(JF\tilde y,\tilde u) +g(JF\tilde x,\tilde u)g(J\tilde y,\tilde z) - 
												g(JF\tilde x,\tilde z)g(J\tilde y,\tilde u) \\
		+2g(F\tilde x,J\tilde y)g(J\tilde z,\tilde u) - 2g(\tilde x,J\tilde y)g(JF\tilde z,\tilde u) \} \ .             
	\end{array}   \leqno (4.1)
$$

Let $X,Y,Z$ be orthogonal tangent vectors at a point $p$ of $M^n$. Then (4.1) and
$ [A_{JX},A_{JY}]=0$ imply

$$
	R(X,Y,Z,X) = \frac \mu 8g(X,X)g(FY,Z) \ .   \leqno (4.2) 
$$

Let \ $X,Y \in T_p(M_1^{n-1}(c))$. Then we find \ $R(X,Y)Z=0$ for any vector \ $Z \in T_pM$,
orthogonal to $X$ and to $Y$. Hence using (4.2) we obtain $g(FY,U)=0$. Consequently
for any \ $Y \in T_pM_1^{n-1}(c)$ \ it follows \ $\pi_1Y=0$ or $\pi_2Y=0$. Suppose now that there
exist nonzero vectors \ $U,V \in T_pM_1^{n-1}(c)$, such that $\pi_1U=0$ and $\pi_2V=0$.
But we must have \ $\pi_1(U+V)=0$ or $\pi_2(U+V) = 0$. Let for example \ $\pi_1(U+V)=0$.
Then $ \pi_1V=0$, which is a contradiction. Consequently we have either $\pi_1=0$ or
$\pi_2=0$ on $T_pM_1^{n-1}(c)$. Hence we obtain easily that $k=n-1$ and \ 
$ M_1^{n-1}(c) \subset \widetilde M_1^{2(n-1)}(\mu)$, $I \subset \widetilde M^2(-\mu)$.
Since $ M_1^{n-1}(c)$ is semiparallel in $\widetilde M_1^{2(n-1)}(\mu)$ and $\mu \ne 0$
it follows that $M_1^{(n-1)}(c)$ is totally geodesic in $\widetilde M_1^{2(n-1)}(\mu)$,
see \cite{OK}, so $\ds c = \frac {\mu}{4}$.

\vspace{0.6in}

\centerline{Somario}

\vspace{0.2in}
{\it In una variet\`a k\"ahleriana si considerano le sottovariet\`a totalmente reali e conforme\-mente
piatte con vettore di curvatuta media parallelo e le sottovariet\`a con seconda forma fondamen\-tale 
semiparallelo.

Sono anche considerate le sottovariet\`a totalmente reali di una variet\`a prodotto 
di variet\`a k\"ahleriane, avente tensore di Bochner nullo.}

\end{document}